\input ./style/arxiv-general.cfg
\documentclass[aop,MSNbibl,dvips]{arximspdf}
\makeatletter
   \@ifpackageloaded{graphicx}{}{\usepackage{graphicx}}
\makeatother
\doi{10.1214/13-AOP860} %kopijuoti is PTS
\volume{43}
\issue{3}
\pubyear{2015}
\firstpage{971}
\lastpage{991}
\makeatletter
\newcommand{\eqref}[1]{(\ref{#1})}
\newcommand{\R}{\mathbb{R}}
\newcommand{\E}{\mathbb{E}}
\newcommand{\curlyS}{\mathcal{S}}
\newcommand{\normal}{\mathcal{N}}
\newcommand{\grad}{\nabla}
\newcommand{\symdiff}{\Delta}
\renewcommand{\Pr}{\operatorname{Pr}}

\newcommand{\tr}{\operatorname{tr}}
\newcommand{\Hess}{H}
\newcommand{\sgn}{\operatorname{sgn}}

\newtheorem{theorem}{Theorem}[section]
\newproclaim{assumption}{Assumption}[section]

\newtheorem{lemma}[theorem]{Lemma}

\newtheorem{proposition}[theorem]{Proposition}
\newproclaim{remark}[theorem]{Remark}
\newproclaim{note}[theorem]{Note}
\newproclaim{definition}[theorem]{Definition}

\makeatother

\begin{document}
\begin{frontmatter}

\title{Robust dimension free isoperimetry in~Gaussian~space\thanksref{T1}}
\runtitle{Robust dimension free Gaussian isoperimetry}

\begin{aug}
\author[A]{\fnms{Elchanan} \snm{Mossel}\ead[label=e1]{mossel@stat.berkeley.edu}\corref{}}
\and
\author[A]{\fnms{Joe} \snm{Neeman}\ead[label=e2]{joeneeman@gmail.com}}
\thankstext{T1}{Supported by Grants NSF (DMS-11-06999) and ONR
(DOD ONR N000141110140).}
\runauthor{E. Mossel and J. Neeman}
\affiliation{University of California, Berkeley}
\address[A]{University of California, Berkeley\\
367 Evans Hall\\
Berkeley, California 94720\\
USA\\
\printead{e1}\\
\phantom{E-mail:\ }\printead*{e2}}
%\author[A]{\fnms{} \snm{}}
%\and
%\author[B]{\fnms{} \snm{}}
%\runauthor{}
%\affiliation{}
%\dedicated{}
%\address[A]{} %adresu isvedimo komanda gale!
%\address[B]{}
\end{aug}

% HISTORY:
\received{\smonth{11} \syear{2012}}
\revised{\smonth{5} \syear{2013}}
%\accepted{\smonth{} \syear{}}

% ABSTRACT
%
\begin{abstract}
We prove the first robust dimension free isoperimetric result for the
standard Gaussian measure $\gamma_n$ and the corresponding boundary
measure $\gamma_n^{+}$ in~$\R^n$.
The main result in the theory of Gaussian isoperimetry (proven in the
1970s by Sudakov and Tsirelson, and independently by Borell) states that
if $\gamma_n(A) = 1/2$ then the surface area of $A$ is bounded by the
surface area of a half-space with the same measure, $\gamma_n^+(A)
\leq(2\pi)^{-1/2}$.
Our results imply in particular that if $A \subset\R^n$ satisfies
$\gamma_n(A) = 1/2$ and
$\gamma_n^+(A) \leq(2\pi)^{-1/2}+\delta$
then there exists a half-space $B \subset\R^n$ such that
$\gamma_n(A \Delta B) \leq C \smash{\log^{-1/2}} (1/\delta)$
for an absolute constant $C$.
Since the Gaussian isoperimetric result was established, only recently
a robust version of the Gaussian isoperimetric result was obtained by
Cianchi et al.,
who showed that $\gamma_n(A \symdiff B) \le C(n) \sqrt\delta$ for
some function $C(n)$ with no effective bounds.
Compared to the results of Cianchi et al.,
our results have optimal (i.e., no) dependence on the dimension,
but worse dependence on $\delta$.
\end{abstract}

% KEYWORDS

\begin{keyword}
\kwd{Noise stability}
\kwd{Gaussian measure}
\kwd{isoperimetric inequalities}
\kwd{majority is stablest}
\end{keyword}
%
%\begin{keyword}
%\kwd[Primary ]{}
%\kwd{}
%\kwd[; secondary ]{}
%\end{keyword}
%\begin{keyword}
%\kwd{}
%\end{keyword}

\end{frontmatter}

%s1 #&#
\section{Introduction}\label{sec1}

Gaussian isoperimetric theory is an extensive and rich theory. It connects
numerous areas of mathematics including probability, geometry~\cite
{MilmanSchechtman86},
concentration
and high dimensional phenomena \cite{ledouxbook}, rearrangement
inequalities \cite{ck01}
and more. For an introduction to Gaussian isoperimetry and its many
applications,
see Ledoux's St.-Flour lecture notes \cite{Ledoux96}.

The main result in this area is that half-spaces minimize the surface
area among all sets with a given Gaussian measure.
This fact, originally proven by Sudakov and Tsirelson \cite{st78}
and independently by Borell \cite{bor75}, now has several
other proofs: Ehrhard \cite{Ehrhard83,Ehrhard84,Ehrhard86}
developed a symmetrization technique, Bakry and
Ledoux \cite{BakryLedoux96} and Ledoux \cite{Ledoux98} used semigroup methods
and Bobkov \cite{bob97} gave a proof based on an isoperimetric
inequality on the discrete cube.

Some of the strongest results in this area deal with extensions
of this basic theorem. For example, Borell \cite{bor85} proved
that half-spaces minimize a more global version of surface area called
Gaussian noise sensitivity. This fact has recently found
applications in areas such as quantitative social choice
and theoretical computer science; see, for example, \cite
{KKMO04,MoOdOl10,Raghavendra08,Mossel12}.

A second direction of extension is characterizing the case
of equality or almost equality. It took about a decade until the
equality case was addressed.
Erhard \cite{Ehrhard86} showed that if a sufficiently nice set
achieves equality then it is a half-space.
It took 15 more years until the same result was proven for general
measurable sets by Carlen and Kerce \cite{ck01}.

Prior to our work, the harder question of almost equality was only
recently addressed by
Cianchi et al. \cite{italian11} who showed that if the Gaussian
boundary of
set $A$ is within $\delta$ of the optimal value then there exists
a half space $B$ such that the Gaussian measure of the symmetric
difference between $A$ and $B$ is at most $C(n) \sqrt\delta$.
Their result gives no bound on the function $C(n)$. Indeed the
techniques of \cite{italian11} are not appropriate for deriving any
bound on $C(n)$.

Our goal in this paper is to establish a robustness result that is
\emph{dimension independent}. Not only such result is more elegant,
it is much in the spirit of Gaussian isoperimetric theory, where the
statement of most results are dimension independent. In particular, the
results of Sudakov and Tsirelson \cite{st78}
and Borell \cite{bor75} are dimension independent: the bound they give
on the size of the boundary in terms of the measure of the set are
dimension independent. Similarly, the results of Borell~\cite{bor85}
are stated in a dimension-free way.
%In particular
%we note that in many of the recent applications of Borell's results
%mentioned above, it is crucial that there is no dependence on the
%dimension. Indeed,
%this is one of the main motivations for the results of the current
%paper.

%s1.1 #&#
\subsection{Gaussian isoperimetry}
The Gaussian isoperimetric inequality was first proved by Sudakov and
Tsirelson \cite{st78},
and independently by Borell \cite{bor75}. It states that in $\R^n$
with the standard Gaussian measure,
the isoperimetric sets are half-spaces. To be more precise, let $\phi
(x) = (2\pi)^{-1/2} e^{-x^2/2}$
be the standard Gaussian density, and define $\Phi(x) = \int_{-\infty
}^x \phi(y) \,dy$.
Let $\gamma_n$ be the standard Gaussian measure on $\R^n$, and define
the boundary measure
$\gamma_n^+$ by
\[
\gamma_n^+(A) = \sup \biggl\{ \int_A (\grad-
x) \cdot v(x) \,d\gamma_n \dvtx v \in\curlyS\bigl(\R^n,
\R ^n\bigr), \bigl|v(x)\bigr| \le1 \mbox{ for all } x \biggr\},
\]
where $\curlyS(\R^n, \R^n)$ is the set of smooth functions
$\R^n \to\R^n$ such that all derivatives vanish at infinity.
(This definition of boundary measure coincides with Minkowski content
and the $(n-1)$-dimensional
Gaussian-weighted Hausdorff measure for sufficiently nice sets \cite{ck01}.)
Then the Gaussian isoperimetric inequality
states that for every measurable $A$,
$\phi(\Phi^{-1}(\gamma_n(A))) \le\gamma_n^+(A)$.
It is not hard to verify that equality is attained if $A$ is an affine
half-space
(i.e., a set of the form $\{x \in\R^n\dvtx x \cdot a \ge b\}$).
In the theory, it is common to define the isoperimetric profile $I =
\phi\circ\Phi^{-1}$, so
that the Gaussian isoperimetric inequality reads
%
%e1 #&#
\begin{equation}
\label{eqisoperimetric} I\bigl(\gamma_n(A)\bigr) \le\gamma_n^+(A).
\end{equation}

%s1.2 #&#
\subsection{Bobkov's inequality}
Bobkov's inequality \cite{bob97} is a functional generalization of the
Gaussian isoperimetric inequality.
The equality case was proved by Carlen and Kerce \cite{ck01}.
Here and for the rest of this article, we will write ``$\E$'' for the
integral with respect
to $\gamma_n$ and $\| \cdot\|$ for the Euclidean norm on $\R^n$.

%th1.1 #&#
\begin{theorem}[(\cite{bob97,ck01})]\label{thmbobkov}
For any smooth function $f\dvtx\R^n \to[0, 1]$ of bounded variation,
%
%e2 #&#
\begin{equation}
\label{eqbobkov} I(\E f) \le\E\sqrt{I^2(f) + \|\grad f\|^2}.
\end{equation}
Equality is attained only if $f(x) = \Phi(a \cdot x + b)$ for some $a
\in\R^n$, $b \in\R$.
\end{theorem}

Using standard approximation techniques \cite{ck01}, one can make
sense of Theorem~\ref{thmbobkov} for functions $f$ that are not smooth.
In particular,
it is possible to take $f$ to be the indicator
function of a set $A$; in that case, $I(f)$ is identically zero
and so~\eqref{eqbobkov} becomes
\[
I\bigl(\gamma_n(A)\bigr) \le\E\|\grad1_A\| =
\gamma_n^+(A)
\]
whenever $A$ is nice enough.
This is just the Gaussian isoperimetric inequality again.
(We will make the above connection rigorous in
Section~\ref{secset}.)
In this limiting case, the nonsmooth equality cases $1_{\{a \cdot x + b
\ge0\}}$
appear. These are easily seen to be limits of the equality cases
$\Phi(a \cdot x + b)$ in Theorem \ref{thmbobkov}.

%s1.3 #&#
\subsection{Robustness}
Our goal in this article is to study the robustness of Theorem~\ref
{thmbobkov}: suppose
that we have a function $f$ which almost achieves equality. Must there
be some
$a$ and $b$ for which $f$ is close to a function of the form $x \mapsto
\Phi(a \cdot x + b)$? For the case of sets---which is perhaps the most
interesting
case---this question was previously studied by Cianchi et al. \cite
{italian11}, who gave a dimension-dependent
estimate.

%th1.2 #&#
\begin{theorem}[(\cite{italian11})]\label{thmitalian}
If $A \subset\R^n$ satisfies $I(\gamma_n(A)) \ge\gamma_n^+(A) -
\delta$, then there is
a half-space $B$ such that
\[
\gamma_n(A \symdiff B) \le C\bigl(n, \gamma_n(A)\bigr)
\sqrt\delta,
\]
where $C(n, r)$ is some function of $n$ and $r$.
\end{theorem}

Due to use of compactness arguments, there are no effective bounds on
the function $C(n,r)$.

Theorem \ref{thmitalian} is sharp in its $\delta$-dependence;
however, the $n$-dependence is certainly
not sharp. Indeed, one often finds things in Gaussian space to be independent
of the dimension.
The isoperimetric inequality itself is an example of this phenomenon, as
the Gaussian isoperimetric profile $I$ does not depend on the dimension.

Note
that the situation is quite different in Euclidean space with the
Lebesgue measure---for which the techniques
used in \cite{italian11} were originally developed---where the
isoperimetric profile $x \mapsto n \omega_n^{1/n} x^{(n-1)/n}$
does depend on $n$.

%s1.4 #&#
\subsection{Our results: Robust and dimension-free}

Dimension-free estimates are satisfying in themselves, but they are also
crucial for certain applications.
As an example, consider Borell's noise stability inequality \cite{bor85}:
take $X, Y \in\R^n$ jointly Gaussian
with $X, Y \sim\normal(0, I)$ and $\E X_i Y_j = \rho\delta_{ij}$.
Then $\Pr(1_A(X) \ne1_A(Y))$ is minimized, over all sets $A$ with prescribed
Gaussian volume, by affine half-spaces. Ledoux showed \cite{ledoux94} that
this generalizes the Gaussian isoperimetric inequality, which is
recovered in the limit as $\rho\to1$. As mentioned above, for
applications of this result it is crucial that is dimension-free.

A robust dimension-free version of Borell's result would immediately
imply a number of important results.
For example, it would show that if a balanced low influence Boolean
function is almost as stable as the majority function, then the
function is close to a weighted majority. Similarly, it would show that
if a balanced low influence function has a Condorcet paradox
probability that is almost as small as that of a majority then it must
be close to a weight majority of a subset of the coordinates.
(Both of the statements above follow from the arguments of \cite{MoOdOl10}.)

The potential applications above further motivate our main result which
is a dimension-independent stability result for Bobkov's inequality in
Gaussian space.
We note, however, that our dependence on $\delta$ is much worse than the
one in Theorem \ref{thmitalian}; improving this dependence is
therefore a natural
open problem.

Our main functional result is the following.

%th1.3 #&#
\begin{theorem}\label{thmmain}
There exists a universal constant $C$ such that the following holds.
Let $f\dvtx\R^n \to\R$ be a smooth function and define
\[
\delta= \E\sqrt{I^2(f) + \|\grad f\|^2} - I(\E f).
\]
There exists a function $g$ of
the form $g(x) = \Phi(a \cdot x + b)$ such that
\[
\E(f - g)^2 \le C \frac{1}{\sqrt{\log(1/\delta)}}.
\]
\end{theorem}

Of course, the most interesting special case of Theorem \ref{thmmain}
is when $f$ is the indicator function of some set. Such an $f$ is not
smooth, of course, but the same arguments that reduced Theorem \ref{thmbobkov}
to the Gaussian isoperimetric inequation can be employed here.
Thus, we obtain a robustness result for the Gaussian
isoperimetric inequality.

%th1.4 #&#
\begin{theorem}\label{thmmain2}
There exists an absolute constant $C$ such that the following holds.
For any measurable set $A \subset\R^n$, let
$\delta= \gamma_n^+(A) - I(\gamma_n(A))$. There exists an affine half-space
$B$ such that
\[
\gamma_n(A \symdiff B) \le C \frac{1}{\sqrt{\log(1/\delta)}}.
\]
\end{theorem}

%s1.5 #&#
\subsection{Proof techniques}
Our approach builds on the work of Carlen and Kerce~\cite{ck01} (which
extends ideas of Ledoux \cite{ledoux94}).
Carlen and Kerce \cite{ck01} write an integral formula [equation (\ref
{eqck}) below] which
bounds $\delta(f) = \E\sqrt{I^2(f) + \|\grad f\|^2} - I(\E f)$ from below.

The ``main term'' in the integral is the Frobenius norm of the Hessian
of $h_t = \Phi^{-1} \circ(P_t f)$,
where $P_t$ is Ornstein--Uhlenbeck semigroup. It is easy to verify that
if $f$ is an indicator of a half-space or
if $f = \Phi(a \cdot x + b)$
then $h_t$ is linear. Our first step in the proof is to utilize a
second-order Poincar\'e inequality which implies that if the
Frobenius norm of the Hessian of $h_t$ is sufficiently small, then
$h_t$ is close to a linear function.

The main effort in our approach is devoted to controlling the
``secondary terms'' in (\ref{eqck}).
This main effort is established in a sequence of analytic results using
the smoothness of the semigroup $P_t$ and involving---among other
techniques---concentration of measure and reverse
hypercontractivity.
Using the approach above, we show that if $\delta= \delta(f)$ is
small then there exists some $t$, not too large, such that
$h_t$ is $\varepsilon(\delta)$-close to a linear function.

The next step of the proof requires applying $P_t^{-1} \circ\Phi$ to
conclude that
$f$ is close to a linear function. There is an obvious obstacle in
this approach: $P_t^{-1}$
is not a bounded operator.

Fortunately, using the smoothness of the original function $f$, or the
fact that we may assume that the original set $A$ has small boundary,
we may deduce a decay in the Hermite expansion of $f$ (or $A$). Thus,
we show that
for the functions under consideration, $P_t^{-1}$ is ``effectively''
bounded which allows us to conclude that $f$ is close to a Gaussian (or
$A$ is close to a half space), proving the result.

%A major part of the analysis is devoted to showing
%that if the integral is small, then there exists a finite $t$ such that
%$\Phi^{-1} \circ(P_t f)$ is close (as a function of $\delta$) to a
%linear function, where $P_t$ is Ornstein-Uhlenbeck semigroup.
%This part uses a number of analytic techniques including delicate
%estimates on derivatives, a second order Poincare inequality and
%the use of reverse-hyper-contraction and concentration.

%The remaining step requires applying $P_t^{-1}$ which in principle
%is not a bounded operator. Using the smoothness
%(or small boundary) of the function $f$, it is possible to truncate
%to the low-degree Hermite expansion and obtain the required result.

%s2 #&#
\section{Semigroup proof of Bobkov's inequality}

Our work begins with Ledoux's short and elementary proof \cite
{ledoux94} of \eqref{eqbobkov}.
The main ingredient of this proof is the Ornstein--Uhlenbeck semigroup:
for $t \ge0$, define the operator $P_t\dvtx L_\infty(\R^n) \to
L_\infty
(\R^n)$ by
\[
(P_t f) (x) = \int_{\R^n} f\bigl(e^{-t}
x + \sqrt{1 - e^{-2t}} y\bigr) \, d\gamma_n(y).
\]
Clearly, $P_0$ is the identity operator and $P_t f$
converges pointwise to $\E f$ as $t \to\infty$.
Consider, therefore, the quantity
%
%e3 #&#
\begin{equation}
\label{eqbobkov-with-semigroup} \E\sqrt{I^2(P_t f) + \|\grad
P_t f\|^2}.
\end{equation}

When $t = 0$, this is exactly the right-hand side of \eqref
{eqbobkov}; as
$t \to\infty$, it approaches the left-hand side of \eqref{eqbobkov}
by the dominated
convergence theorem and the boundedness of $f$.
To prove \eqref{eqbobkov}, Ledoux
differentiated \eqref{eqbobkov-with-semigroup} with respect to $t$
and showed
that the derivative is nonpositive. Thus, a potentially difficult
inequality turns
into a calculus problem.

Actually, Ledoux only explicitly
differentiated \eqref{eqbobkov-with-semigroup} in the one-dimensional case.
The $n$-dimensional case of \eqref{eqbobkov-with-semigroup} was computed
by Carlen and Kerce \cite{ck01} in their work on the equality case.
Our robustness
result is based on their calculations, which we will summarize as a lemma.

%le2.1 #&#
\begin{lemma}[(\cite{ck01})]\label{lemck}
Let $f$ be a smooth function $f\dvtx\R^n \to[0, 1]$. Define
$h_t = \Phi^{-1} \circ(P_t f)$ and
\[
\delta(f) = \E\sqrt{I^2(f) + \|\grad f\|^2} - I(\E f).
\]
Then
%
%e4 #&#
\begin{equation}
\label{eqck} \delta(f) \ge\int_0^\infty\E
\frac{\phi(h_t) \|H(h_t)\|_F^2}{(1 + \|\grad
h_t\|^2)^{3/2}} \,dt,
\end{equation}
where $H(h_t)$ is the Hessian matrix of $h_t$
and $\|\cdot\|_F$ denotes the Frobenius norm $\|A\|_F^2 = \tr(A^T A)$.
\end{lemma}

From now on, $\delta(f)$ will be defined as it was in Lemma \ref{lemck}.
Where $f$ is clear from the context, we will only write $\delta$.

The equality case in Theorem \ref{thmbobkov} follows fairly easily from
Lemma \ref{lemck}: if \mbox{$\delta= 0$}, then $H(h_t)$ must be zero for
all $t > 0$,
which implies that $h_t$ is a linear function for all $t > 0$. A
straightforward limiting
argument shows that one can take $t$ to zero, and the result follows.

Our proof of Theorem \ref{thmmain} works by finding a lower bound for
the right-hand side of \eqref{eqck}.
First, we replace the integral over $[0, \infty)$
by an integral over $[C, C+1]$, where $C$ is a large enough constant.
For some $t \in[C, C+1]$, we find an affine function $h^*$ such that
$h_t$ is close to $h^*$. In particular, this means that $f_t$
is close to $\Phi\circ h^*$. This part of the argument will be carried
out in Section~\ref{seclarge-t}.
The second part of the argument, carried out in Section~\ref{secsmall-t},
shows that $f$ must be close to $P_{-t} (\Phi\circ h^*)$.

%s3 #&#
\section{Approximation for large $t$}\label{seclarge-t}

This section is devoted to the proof of Proposition \ref
{proplarge-t}, which
shows that $h_t$ can be approximated by an affine function for some
sufficiently large $t$.

%pr3.1 #&#
\begin{proposition}\label{proplarge-t}
There is a universal constant $C > 0$ such that
for any measurable $f\dvtx\R^n \to[0, 1]$ there exists $t \in[C, C+1]$
such that
\[
\E\bigl(h_t(X) - \E h_t - X \cdot\E\grad
h_t\bigr)^2 \le C \frac{\delta^{1/4}(f)}{m(f)^{5/4}},
\]
where $h_t = \Phi^{-1} \circ(P_t f)$ and $m(f) = (\E f) (1 - \E f)$.
\end{proposition}

%no3.2 #&#
\begin{note}
In this section and the next, we will not be concerned with the value of
universal constants; hence, the letters $C$ and $c$ will denote universal
constants, whose values may change from line to line. We will use $C$ to
denote constants that must be sufficiently large, while $c$ denotes
constants that must be sufficiently small.
\end{note}

\textit{Notation}.
From now on, we fix the notation
\[
f_t:= P_t f, \qquad h_t:= \Phi^{-1}
\circ f_t,\qquad k_t = e^{-t} /
\sqrt{1-e^{-2t}}.
\]

%s3.1 #&#
\subsection{A second-order Poincar\'e inequality}

In proving the equality cases of Bobkov's inequality, Carlen and Kerce
used the fact that if
$\|\Hess(h_t)\|_F^2$ vanishes then $h_t$ must be a linear function.
The first step
toward Proposition \ref{proplarge-t} is a quantitative version of
this observation.

%le3.3 #&#
\begin{lemma}\label{lempoincare-hessian}
For any smooth function $h\dvtx\R^n \to\R$,
\[
\E\bigl(h(X) - \E h - X\cdot\E\grad h\bigr)^2 \le\E\bigl\|\Hess(h)
\bigr\|_F^2.
\]
\end{lemma}

Since $x \mapsto\E h + x \cdot\E\grad h$ is a linear function,
Lemma \ref{lempoincare-hessian} implies that if $\E\|\Hess(h)\|_F^2$
is small then $h$ is close to linear.
This puts us on our way toward the proof of Proposition \ref{proplarge-t}.
Indeed, if we could remove the $\phi(h_t) (1 + |\grad h_t|^2)^{-3/2}$
term from
the right-hand side of \eqref{eqck}, we would be done already.
The removal of this nuisance term is the topic of the next section.

\begin{pf*}{Proof of Lemma \ref{lempoincare-hessian}}
Recall Poincar\'e's inequality \cite{Ledoux00}
%
%e5 #&#
\begin{equation}
\label{eqpoincare} \E h^2 - (\E h)^2 \le\E\|\grad h
\|^2.
\end{equation}
If we apply \eqref{eqpoincare} to the partial derivatives of $h$, we obtain
\[
\E \biggl(\frac{\partial h}{\partial x_i} \biggr)^2 - \biggl(\E\frac
{\partial h}{\partial x_i}
\biggr)^2 \le\E \biggl(\frac{\partial h}{\partial x_i} \biggr)^2 \le\sum
_{j=1}^n \E \biggl(\frac{\partial^2 h}{\partial x_i\,
\partial x_j}
\biggr)^2.
\]
Summing over $i$ yields $\E\|\grad h\|^2 - \|\E\grad h\|^2
\le\E\|\Hess(h)\|_F^2$.
Combining this with~\eqref{eqpoincare}, we have
\[
\E(h - \E h - X\cdot\E\grad h)^2 = \E h^2 - (\E
h)^2 - \|\E\grad h\|^2 \le\E\bigl\|H(h)\bigr\|_F^2,
\]
where the first equality follows because integration by parts implies
that  $\E X h(X) = \E\grad h$; hence, the orthogonal projection of $h$
onto the span
of linear functions is $X\cdot\E\grad h$.
\end{pf*}

%s3.2 #&#
\subsection{First derivative bounds}

We have shown how to use the $\E\|H(h_t)\|_F^2$ term on the right-hand
side of \eqref{eqck}. In this section, we discuss the
$(1 + \|\grad h_t\|^2)^{-3/2}$ term. A result by Bakry and
Ledoux \cite{BakryLedoux96} shows that this term may be bounded pointwise
from below.

%th3.4 #&#
\begin{theorem}[(\cite{BakryLedoux96})]\label{thmgrad-bound}
For any measurable function $f\dvtx\R^n \to[0, 1]$ and any $t > 0$,
\[
\|\grad f_t\| \le k_t I(f_t)
\]
pointwise,
where $k_t = e^{-t} / \sqrt{1 - e^{-2t}}$ and $f_t = P_t f$.
Equivalently,
\[
\|\grad h_t\| \le k_t
\]
pointwise, where $h_t = \Phi^{-1} \circ f_t$.
\end{theorem}

Note that the second inequality is equivalent to the first by the chain
rule, since $\frac{d }{d x} \Phi^{-1}(x) = \frac{1}{I(x)}$.
Note also that since $I(x) \sim x \sqrt{2\log(1/x)}$ as $x \to0$,\vspace*{1pt}
Theorem~\ref{thmgrad-bound} follows, up to a constant factor, from
the reverse log-Sobolev inequality~\cite{Ledoux00}
%
%e6 #&#
\begin{equation}
\label{eqreverse-log-sob} \frac{\|\grad f_t\|^2}{2 k_t^2 f_t} \le P_t (f \log f) -
f_t \log f_t.
\end{equation}
However, the sharp constant in Theorem \ref{thmgrad-bound} will
be useful in Section~\ref{secsmall-t}.

%s3.3 #&#
\subsection{Reverse-hypercontractivity and reverse-H\"older}

Recall our current task: a~lower bound on \eqref{eqck} for large
$t$. We have already shown that $\|\grad h_t\|$ must be small
for large $t$; our goal for this section is to find a lower bound
on $\E\phi(h_t) \|H(h_t)\|_F^2$.

%pr3.5 #&#
\begin{proposition}\label{proppost-reverse}
There exists a constant $C$ such that if $t \ge C$ then
\[
\E \bigl( \phi(h_t) \bigl\|H(h_t)\bigr\|_F^2
\bigr) \ge \tfrac{1}4 \bigl( I(\E f) \bigr)^2 \bigl( \E
\bigl\|H(h_t)\bigr\|_F \bigr)^2.
\]
\end{proposition}

For this, we will use two inequalities: Borell's
reverse-hypercontractive inequality and the reverse-H\"older
inequality. The reverse-H\"older inequality is classical: for any
$p < 1$ and any positive functions $f$ and $g$,
%
%e7 #&#
\begin{equation}
\label{eqrev-holder} \E fg \ge \bigl(\E f^p \bigr)^{1/p} \bigl(\E
g^{p/(p-1)} \bigr)^{(p-1)/p}.
\end{equation}
The reverse-hypercontractive inequality was proved by Borell \cite{bor82}:
for a positive function $f$ and any $p < 1$, $t > 0$,
%
%e8 #&#
\begin{equation}
\label{eqrev-hyper} \bigl(\E(P_t f)^p \bigr)^{1/p}
\ge \bigl(\E f^q \bigr)^{1/q},
\end{equation}
where $q = 1 + e^{-2t}(p-1)$.

\begin{pf*}{Proof of Proposition \ref{proppost-reverse}}
Let $g_t = \|H(h_t)\|_F^2$ and
apply \eqref{eqrev-holder} with $p = \frac{1}{2}$:
%
%e9 #&#
\begin{equation}
\label{eqapply-rev-holder} \E \bigl(\phi(h_t) \bigl\|H(h_t)
\bigr\|_F^2 \bigr) = \E \bigl(I(f_t)
g_t \bigr) \ge \bigl(\E\bigl(I(f_t)
\bigr)^{-1} \bigr)^{-1} (\E\sqrt g_t
)^2.
\end{equation}
Now, $I$ is a concave function, and so
$I(P_t f) \ge P_{t/2} I(P_{t/2} f)$. Applying \eqref{eqrev-hyper}
with $p = -1$ gives
%
%e10 #&#
\begin{equation}
\label{eqapply-rev-hyper} \bigl(\E\bigl(I(f_t)\bigr)^{-1}
\bigr)^{-1} \ge \bigl(\E\bigl(P_{t/2} I(f_{t/2})
\bigr)^{-1} \bigr)^{-1} \ge \bigl(\E\bigl(I(f_{t/2})
\bigr)^q \bigr)^{1/q},
\end{equation}
where $q = 1 - 2e^{-t}$.
If $t \ge2$, then $q \ge1/2$.
Hence, we can combine \eqref{eqapply-rev-holder}
with \eqref{eqapply-rev-hyper} to obtain
\[
\E \bigl(I(f_t) g_t \bigr) \ge \bigl(\E
\sqrt{I(f_{t/2})} \bigr)^2 (\E\sqrt g_t
)^2.
\]

It remains to show that
%
%e11 #&#
\begin{equation}
\E\sqrt{I(f_{t/2})} \ge\tfrac{1}{2} I(\E f). \label{eqpost-reverse-goal}
\end{equation}
Applying \eqref{eqbobkov} to $f_{t/2}$, we have
%
%e12 #&#
\begin{equation}
\label{eqbobkov-at-t} I(\E f) = I(\E f_{t/2}) \le\E I(f_{t/2}) + \E\|
\grad f_{t/2}\|,
\end{equation}
while Theorem \ref{thmgrad-bound} gives
\[
\E\|\grad f_{t/2} \| \le k_{t/2} \E I (f_{t/2}).
\]
For $t$ sufficiently large, $k_{t/2} \le1$ and so
$\E\|\grad f_{t/2}\| \le\E I(f_{t/2})$;
by \eqref{eqbobkov-at-t}, $I(\E f) \le2 \E I(f_{t/2})$
for large enough $t$.
Now, $I$ is bounded above by $(2\pi)^{-1/2} \le1$, and so
$\E\sqrt{I(f_{t/2})} \ge\E I(f_{t/2}) \ge\frac{1}2 I(\E f)$,
which proves \eqref{eqpost-reverse-goal} and the proposition.
\end{pf*}

%s3.4 #&#
\subsection{Second-derivative estimates}

There is one more ingredient in the proof of Proposition \ref{proplarge-t}:
an upper bound on second derivatives of $h_t$. To see why such a bound
is useful, note that
Lemma \ref{lempoincare-hessian} gives a lower bound on
$\E\|H(h_t)\|_F^2$, but Proposition \ref{proppost-reverse}
contains $\E\|H(h_t)\|_F$. To combine these two results, we
must therefore bound the first moment
of $\|H(h_t)\|_F$ from below in terms of the second moment.
This can be done by H\"older's inequality, as long as we can bound
higher moments of $\|H(h_t)\|_F$ from above. Such a bound
is the goal of this section.

The main bound of this section is the following proposition.

%pr3.6 #&#
\begin{proposition}\label{prophess-bound-lp}
There is a constant $C$ such that for all $t > C$,
\[
\bigl(\E\bigl\|H(h_t)\bigr\|_F^3
\bigr)^{1/3} \le\sqrt{\log\frac{1}{m(f)}},
\]
where $m(f) = \E f (1 - \E f)$.
\end{proposition}

Proposition \ref{prophess-bound-lp} essentially follows by
integrating a pointwise bound on
$\|H(h_t(x))\|_F$:

%le3.7 #&#
\begin{lemma}\label{lemhess-bound}
There is a constant $C$ such that for all $x \in\R^n$
and $t > 0$,
\[
\bigl\|H\bigl(h_t(x)\bigr)\bigr\|_F \le C k_t^2
\sqrt{\log\frac{1}{f_t(x)(1 - f_t(x))}},
\]
where $k_t = e^{-t}/\sqrt{1-e^{-2t}}$.
\end{lemma}

And we will also need to relate the median of $f_t$ with its mean.

%le3.8 #&#
\begin{lemma}\label{lemmedian-mean}
If $M_t$ is a median of $f_t$, then
\[
\E f (1 - \E f) \le2 M_t^{ ({1}/{(1 + k_t)} )^2}.
\]
\end{lemma}

Before we prove either of these lemmas, we will show how they
imply Proposition \ref{prophess-bound-lp}.

\begin{pf*}{Proof of Proposition \ref{prophess-bound-lp}}
Let $g_t = \sqrt{\log(1/f_t)}$ (where $f_t = P_t f$);
note that $\grad g_t = - \frac{\grad f_t}{2 f_t \sqrt{\log(1/f_t)}}$,
and so the reverse log-Sobolev inequality \eqref{eqreverse-log-sob} implies
that $g_t$ is Lipschitz with constant $k_t/\sqrt2$.
Take $N$ to be a median of $g_t$.
By Gaussian concentration for Lipschitz functions,
\[
\E|g_t - N|^p \le\int_0^\infty
\gamma_n\bigl\{|g_t - N|^p \ge x\bigr\} \,dx
\le\int_0^\infty e^{-{x^{2/p}}/{(2 k_t^2)}} \,dx.
\]
After the change of variables $y^2 = \frac{x^{2/p}}{k_t^2}$, the
right-hand side is
just $p k_t^p \E|Y|^{p-1}$, where $Y$ is a standard Gaussian variable.
Since $\E|Y|^{p-1} \le(Cp)^{p/2}$, it follows that
\[
\bigl(\E|g_t - N|^p\bigr)^{1/p} \le\bigl( (C
p)^{{p}/2 + 1} k_t^p\bigr)^{1/p} \le C
k_t \sqrt p.
\]
Then, by the triangle inequality,
%
%e13 #&#
\begin{equation}
\label{eqtriangle} \bigl(\E|g_t|^p\bigr)^{1/p} \le
N + C k_t \sqrt p.
\end{equation}

By Lemma \ref{lemhess-bound}, $\|H(h_t)\|_F \le C k_t^2 g_t$
pointwise whenever $f_t \le\frac{1}{2}$.
Thus,
%
%e14 #&#
\begin{equation}
\label{eqhess-bound-1} \bigl(\E\bigl\|H(h_t)\bigr\|_F^p
1_{\{f_t \le{1}/{2}\}} \bigr)^{1/p} \le C k_t^2
\bigl(\E|g_t|^p \bigr)^{1/p} \le C
k_t^2 (N + k_t \sqrt p),
\end{equation}
where the second inequality follows from \eqref{eqtriangle}.

To relate $N$ to $\E f$, simply
note that $M = e^{-N^2}$ is a median for $f_t$; hence,
Lemma \ref{lemmedian-mean} implies that
$N \le C (1 + k_t) \sqrt{ \log(1/m(f)) }$,
where $m(f) = \E f (1 - \E f)$.
Plugging these bounds into \eqref{eqhess-bound-1},
%
%e15 #&#
\begin{equation}
\label{eqhalf-done} \bigl(\E\bigl\|H(h_t)\bigr\|_F^p
1_{\{f_t \le{1}/{2}\}} \bigr)^{1/p} \le C k_t^2
\biggl( ( 1 + k_t) \sqrt{\log\frac{1}{m(f)}} + k_t
\sqrt p \biggr).
\end{equation}

By the same argument with $f_t$ and $1-f_t$ exchanged, we have
%
%e16 #&#
\begin{equation}
\label{eqother-half} \bigl(\E\bigl\|H(h_t)\bigr\|_F^p
1_{\{f_t \ge{1}/{2}\}} \bigr)^{1/p} \le C k_t^2
\biggl( ( 1 + k_t) \sqrt{\log\frac{1}{m(f)}} + k_t
\sqrt p \biggr).
\end{equation}
Combining \eqref{eqhalf-done}
and \eqref{eqother-half} with the triangle inequality
\begin{eqnarray*}
\bigl(\E\bigl\|H(h_t)\bigr\|_F^p
\bigr)^{1/p} &\le& \bigl(\E\bigl\|H(h_t)\bigr\|_F^p
1_{\{f_t \ge{1}/{2}\}} \bigr)^{1/p} + \bigl(\E\bigl\|H(h_t)
\bigr\|_F^p 1_{\{f_t \le{1}/{2}\}} \bigr)^{1/p}
\\
&\le& C k_t^2 \biggl( ( 1 + k_t) \sqrt{\log
\frac{1}{m(f)}} + k_t \sqrt p \biggr).
\end{eqnarray*}
Setting $p = 3$ and taking $t$ sufficiently large completes the proof.
\end{pf*}

The fact that $\sqrt{\log(1/f_t)}$ is Lipschitz was noticed
by Hino \cite{Hino02}, and was also used recently by
Ledoux \cite{Ledoux11}. This fact, which was important
in the preceding proof, will also be crucial in the proof of
Lemma \ref{lemmedian-mean}:

\begin{pf*}{Proof of Lemma \ref{lemmedian-mean}}
Let $g_t = \sqrt{\log(1/f_t)}$; take $N$ to be a median
of $g_t$ and let $M = e^{-N^2}$, so that $M$ is a median
of $f_t$. For any $\alpha< 1$,
\[
\Pr\bigl(f_t \ge M^{\alpha^2}\bigr) = \Pr\bigl(g_t
\le\alpha\sqrt{\log(1/M)}\bigr) = \Pr(g_t \le\alpha N).
\]
Recall that $g_t$ is $\frac{1}{\sqrt2} k_t$-Lipschitz.
Thus,
\[
\Pr\bigl(f_t \ge M^{\alpha^2}\bigr) = \Pr(g_t \le
\alpha N) \le\exp \biggl(-\frac{(1-\alpha)^2 N^2}{k_t^2} \biggr) =
M^{{(1 - \alpha)^2}/{k_t^2}}.
\]
Setting $\alpha= \frac{1}{1 + k_t}$, we have
$\frac{(1-\alpha)^2}{k_t^2} = \alpha^2$. Thus,
$\Pr(f_t \ge M^{\alpha^2}) \le M^{\alpha^2}$. Since $f_t \le1$,
Markov's inequality implies that $\E f_t \le2 M^{\alpha^2}$.
\end{pf*}

For the rest of the section, we will devote ourselves to proving
Lemma \ref{lemhess-bound}, which we will do very explicitly.
The proof of Lemma \ref{lemhess-bound} begins with the formula
%
%e17 #&#
\begin{equation}
\label{eq2-diff-h} \frac{\partial^2 h_t}{\partial x_i\, \partial x_j} = \frac{1}{I(f_t(x))}
\frac{\partial^2 f_t}{\partial x_i\, \partial x_j} +
\frac{\Phi^{-1}(f_t(x))}{I^2(f_t(x))} \frac{\partial f_t}{\partial
x_i} \frac{\partial f_t}{\partial x_j}.
\end{equation}
We will bound the two terms on the right-hand side in two different
lemmas. But first, we quote a result on the moments of
a order-2 Gaussian chaos. To obtain Theorem \ref{thmchaos}
from the result stated in \cite{hw71}, simply note that
the operator norm is bounded by the Frobenius norm.

%th3.9 #&#
\begin{theorem}[(\cite{hw71})]\label{thmchaos}
For any matrix $A$ and any $1 \le p < \infty$,
if $Y$ is a standard Gaussian vector in $\R^n$ then
\[
\bigl(\E\bigl|Y^T A Y - \tr(A)\bigr|^p \bigr)^{1/p} \le C
p \|A\|_F.
\]
\end{theorem}

Theorem \ref{thmchaos} will be used to bound the
first term of \eqref{eq2-diff-h}.

%le3.10 #&#
\begin{lemma}\label{lem2-diff-terms}
For any matrix $A = (a_{ij})$,
\[
\sum_{ij} a_{ij} \frac{\partial^2 f_t}{\partial x_i\, \partial x_j} \le C
k_t^2 \|A\|_F f_t \log
\frac{1}{f_t}.
\]
\end{lemma}

\begin{pf}
We write out derivatives of $f_t$
as integrals: if $i \ne j$ then
with the change of variables $y = \frac{z - e^{-t} x}{\sqrt{1 - e^{-2t}}}$,
\begin{eqnarray*}
\frac{\partial^2 f_t}{\partial x_i\, \partial x_j} &=& \int_{\R^n} f(z)
\frac{\partial^2}{\partial x_i\, \partial x_j} \phi
\biggl(\frac{z -e^{-t} x}{\sqrt{1 - e^{-2t}}} \biggr) \frac
{dz}{\sqrt{1 - e^{-2t}}}
\\
&=& k_t^2 \int_{\R^n} f(z)
y_i y_j \phi \biggl(\frac{z -e^{-t} x}{\sqrt{1 - e^{-2t}}} \biggr)
\frac
{dz}{\sqrt{1 - e^{-2t}}}
\\
&=& k_t^2 \int_{\R^n} f
\bigl(e^{-t} x + \sqrt{1-e^{-2t}} y\bigr) y_i
y_j \phi (y) \,dy,
\end{eqnarray*}
while if $i = j$ then a similar computation gives
\[
\frac{\partial^2 f_t}{\partial x_i\, \partial x_j} = k_t^2 \int_{\R^n}
f\bigl(e^{-t} x + \sqrt{1-e^{-2t}} y\bigr)
\bigl(y_i^2 - 1\bigr) \phi(y) \,dy.
\]
Applying H\"older's inequality,
\begin{eqnarray*}
\sum_{ij} a_{ij} \frac{\partial^2 f_t}{\partial x_i\, \partial x_j} &=&
k_t^2 \int_{\R^n} f
\bigl(e^{-t} x + \sqrt{1-e^{-2t}} y\bigr) \bigl(y^T
A y - \tr (A)\bigr)\phi(y) \,dy
\\
&\le& k_t^2 \bigl(P_t f^p
\bigr)^{1/p} \bigl(\E\bigl|Y^T A Y - \tr(A)\bigr|^q
\bigr)^{1/q}
\\
&\le& k_t^2 f_t^{1/p} \bigl(
\E\bigl|Y^T A Y - \tr(A)\bigr|^q\bigr)^{1/q},
\end{eqnarray*}
where $\frac{1}{p} + \frac{1}{q} = 1$
and $Y$ is distributed as a standard Gaussian variable, and the
last line follows because $0 \le f \le1$ and so $f^p \le f$.
By Theorem \ref{thmchaos},
\[
\sum_{ij} a_{ij} \frac{\partial^2 f_t}{\partial x_i\, \partial x_j} \le C
k_t^2 q f_t^{1/p} \|A
\|_F.
\]
Choosing $1/p = 1 - \frac{1}{\log(1/f_t)}$, we have $f_t^{1/p} = e
f_t$ and
$q = \log\frac{1}{f_t}$, proving the claim.
\end{pf}

Putting Lemma \ref{lem2-diff-terms} and Theorem \ref{thmgrad-bound}
together, we arrive at a proof of Lemma~\ref{lemhess-bound}.

\begin{pf*}{Proof of Lemma \ref{lemhess-bound}}
Note that
\[
\bigl\|H(h_t)\bigr\|_F = \sup_{\|A\|_F = 1} \sum
_{ij} a_{ij} \frac{\partial
^2 h_t}{\partial x_i\, \partial x_j}.
\]
For any fixed matrix $A$, with $\|A\|_F = 1$, \eqref{eq2-diff-h}
implies that
\[
\sum_{ij} a_{ij} \frac{\partial^2 h_t}{\partial x_i\, \partial x_j} =
\frac{1}{I(f_t)} \sum_{ij} a_{ij}
\frac{\partial^2 f_t}{\partial x_i\, \partial x_j} + \frac{\Phi^{-1}(f_t)}{I^2(f_t)} (\grad f_t)^T A
\grad f_t.
\]
Now, the reverse log-Sobolev inequality \eqref{eqreverse-log-sob}
implies that
$(\grad f_t)^T A \grad f_t \le\|A\|_F \times\break \|\grad f_t\|^2 \le
2k_t^2 \|A\|_F f_t^2 \log\frac{1}{f_t}$.
This, together with Lemma \ref{lem2-diff-terms} and the fact that
$\|A\|_F = 1$, implies
\[
\sum_{ij} a_{ij} \frac{\partial^2 h_t}{\partial x_i\, \partial x_j} \le C
k_t^2 \biggl( \frac{f_t \log(1/f_t)}{I(f_t)} + \frac{\Phi^{-1}(f_t)}{I^2(f_t)}
f_t^2 \log\frac{1}{f_t} \biggr).
\]
Since the right-hand side is independent of $A$, we can take the
supremum over
$\{A\dvtx\|A\|_F = 1\}$, giving
%
%e18 #&#
\begin{equation}
\label{eqhess-bound1} \bigl\|H(h_t)\bigr\|_F \le C k_t^2
\biggl( \frac{f_t \log(1/f_t)}{I(f_t)} + \frac{\Phi^{-1}(f_t)}{I^2(f_t)} f_t^2
\log\frac{1}{f_t} \biggr).
\end{equation}

Next, we claim that for any $0 < a \le\frac{1}{2}$,
%
%e19 #&#
%e20 #&#
\begin{eqnarray}
a \sqrt{\log(1/a)} &\le& CI(a),\label{eqI-bound}
\\
a^2 \sqrt{\log(1/a)} &\le& C \biggl|\frac{I^2(a)}{\Phi^{-1}(a)} \biggr|. \label{eqI^2Phi-bound}
\end{eqnarray}
Now, \eqref{eqI-bound} follows from the well-known fact
(see, e.g., \cite{BakryLedoux96}) that $I(a) \asymp\break a \sqrt{\log(1/a)}$
as $a \to0$\vspace*{-1pt} [where the notation $f(a) \asymp g(a)$
means that $0 < \liminf\frac{f(a)}{g(a)} \le
\limsup\frac{f(a)}{g(a)} < \infty$]. To
show \eqref{eqI^2Phi-bound}, set $g(a) = -\frac{I^2(a)}{\Phi
^{-1}(a)}$; we then
compute
\[
g'(a) = 2 I(a) + \frac{I(a)}{(\Phi^{-1})^2(a)} \asymp I(a) \asymp a \sqrt{
\log\frac{1}{a}}
\]
as $a \to0$. Since
\[
\frac{d }{d a} a^2 \sqrt{\log\frac{1}{a}} = a \sqrt{\log
\frac{1}{a}} - \frac{a}{2 \sqrt{\log(1/a)}} \asymp a \sqrt{\log\frac{1}{a}},
\]
it follows that $g(a) \asymp a^2 \sqrt{\log(1/a)}$, proving \eqref
{eqI^2Phi-bound}.

Suppose that $f_t \le\frac{1}{2}$. Applying \eqref{eqI-bound}
and \eqref{eqI^2Phi-bound} to \eqref{eqhess-bound1} with $a = f_t$
completes the proof in this case. If $f_t > \frac{1}{2}$, we apply
the same argument to $1 - f_t$.
\end{pf*}

%s3.5 #&#
\subsection{Proof of Proposition \texorpdfstring{\protect\ref{proplarge-t}}{3.1}}

With all of the ingredients laid out, the proof of Proposition \ref
{proplarge-t}
follows easily.

\begin{pf*}{Proof of Proposition \ref{proplarge-t}}
Suppose $C$ is a large enough universal constant so that for all $t \ge C$:
\begin{itemize}
\item$\|\grad h_t\| \le1$ (by Theorem \ref{thmgrad-bound}).
\item$\E (\phi(h_t) \|H(h_t)\|_F^2 )
\ge\frac{1}4 I^2(\E f)  (\E\|H(h_t)\|_F )^2$
(by Proposition \ref{proppost-reverse}).
\item$(\E\|H(h_t)\|_F^3)^{1/3} \le\sqrt{\log(1/(m(f)))}$ (by
Proposition \ref{prophess-bound-lp}),
where $m(f) = \E f (1 - \E f)$.
\end{itemize}
By the first two bullet points, for any $t \ge C$
%
%e21 #&#
\begin{eqnarray}
\label{eqintermediate} \E\frac{\phi(h_t)
 \|H(h_t)\|_F^2}{(1 + \|\grad h_t\|^2)^{3/2}} & \ge&2^{-3/2} \E \bigl(
\phi(h_t) \bigl\|H(h_t)\bigr\|_F^2 \bigr)
\nonumber
\\[-8pt]
\\[-8pt]
\nonumber
& \ge&2^{-7/2} I^2(\E f) \bigl(\E\bigl\|H(h_t)
\bigr\|_F \bigr)^2.
\end{eqnarray}
Now, H\"older's inequality implies that
\[
\E X^2 = \E X^{1/2} X^{3/2} \le (\E X
)^{1/2} \bigl(\E X^3 \bigr)^{1/2}
\]
for any nonnegative random variable $X$. Applying this with
$X = \|H(h_t)\|_F$, the third bullet point above implies that
\begin{eqnarray*}
\E \bigl(\bigl\|H(h_t)\bigr\|_F^2 \bigr) &\le& \bigl(
\E\bigl\|H(h_t)\bigr\|_F \bigr)^{1/2} \bigl(\E
\bigl\|H(h_t)\bigr\|_F^3 \bigr)^{1/2}
\\
&\le& \bigl(\E\bigl\|H(h_t)\bigr\|_F \bigr)^{1/2}
\log^{3/4} \frac{1}{m(f)}.
\end{eqnarray*}
Plugging this into \eqref{eqintermediate},
%
%e22 #&#
\begin{equation}
\label{eqintermediate-2} \E\frac{\phi(h_t)
\|H(h_t)\|_F^2}{(1 + \|\grad h_t\|^2)^{3/2}} \ge2^{-7/2} I^2(\E f)
\bigl(\E \bigl(\bigl\|H(h_t)\bigr\|_F^2 \bigr)
\bigr)^4 \log^{-3} \frac{1}{m(f)}.
\end{equation}
Now, \eqref{eqintermediate-2} holds for any $t \ge C$. In particular,
if we choose $t^* \in[C, C + 1]$ to minimize $\E\|H(h_{t^*})\|_F^2$,
then by Lemma \ref{lemck} and \eqref{eqintermediate-2},
%
%e23 #&#
\begin{eqnarray}
\label{eqintermediate-3} \delta(f) & \ge&\int_0^{\infty} \E
\frac{\phi(h_t) \|H(h_t)\|_F^2}{(1 + \|
\grad h_t\|^2)^{3/2}} \,dt
\nonumber
\\
& \ge&\int_C^{C+1} \E\frac{\phi(h_t) \|H(h_t)\|_F^2}{(1 + \|\grad
h_t\|^2)^{3/2}} \,dt
\nonumber
\\[-8pt]
\\[-8pt]
\nonumber
& \ge&2^{-7/2} \int_{C}^{C+1}
I^2(\E f) \bigl(\E\bigl\|H(h_t)\bigr\|_F^2
\bigr)^4 \log^{-3} \frac{1}{m(f)} \,dt
\nonumber
\\
& \ge&2^{-7/2} \frac{I^2(\E f)}{\log^{3} (1/m(f))} \bigl(\E\bigl\| H(h_{t^*})
\bigr\|_F^2 \bigr)^4.\nonumber
\end{eqnarray}
Recall that $I^2(x) \asymp x^2 \log(1/x)$ as $x\to0$.
Hence, as $x \to0$,
\[
\frac{I^2 (x)}{\log^{3} (1/(x(1-x)))} \asymp\frac{x^4}{\log(1/x)} \ge x^5.
\]
By replacing $x$ with $1-x$ [note that $I(x) = I(1-x)$] and repeating the
argument, we see that there is some universal constant $c > 0$ such that
$I^2(x) \log^{-3} (1/(x(1-x))) \ge c (x(1-x))^5$ for all $x \in[0, 1]$.
Applying this to \eqref{eqintermediate-3} with $x = \E f$, we have
%
%e24 #&#
\begin{equation}
\label{eqdelta-lower-bound} \delta(f) \ge c m(f)^5 \bigl(\E\bigl\|H(h_{t^*})
\bigr\|_F^2 \bigr)^4.
\end{equation}
Finally, by Lemma \ref{lempoincare-hessian},
\[
\E\bigl(h_{t^*}(X) - \E h_t^* - X\cdot\E\grad
h_t^*\bigr)^2 \le\E\bigl\| H(h_{t^*})
\bigr\|_F^2 \le C \frac{\delta^{1/4}(f)}{(\E f)^{5/4}},
\]
where the last inequality follows from \eqref{eqdelta-lower-bound}.
\end{pf*}

%s4 #&#
\section{Approximation for small $t$}\label{secsmall-t}
Recall that $f_t = P_t f$ and $h_t = \Phi^{-1} \circ f_t$.
Proposition \ref{proplarge-t} shows that if $f$
achieves almost-equality in \eqref{eqbobkov} then
$h_t$---for some $t$ not too large---can be well
approximated by a linear function. Since $\Phi$ is a contraction,
this implies that $f_t$ may be well approximated by a function of the form
$\Phi(a\cdot x + b)$. The goal of this section is to
complete the proof of Theorem \ref{thmmain} by
showing that $f$ itself can be approximated by a function
of the same form.
This will be accomplished mainly with spectral techniques,
by expanding $f$ in the Hermite basis.

Let $g_t(x) = \Phi(\E h_t + x \cdot\E\grad h_t)$, so that
Proposition \ref{proplarge-t} implies that
$\E(f_t - g_t)^2 \le C \delta^{1/4}(f) m^{-5/4}(f)$.
By directly computing $P_t$ applied to the indicator of a half-space,
one may check the following lemma (which also appeared implicitly
in \cite{ck01}):

%le4.1 #&#
\begin{lemma}\label{lempullback}
If $\|a\| \le k_t$, then the function $g(x) = \Phi(a\cdot x + b)$
is in the range of $P_t$. Moreover, if $\|a\| = k_t$ then $P_t^{-1} g$
is the indicator function of a half-space, while if $\|a\| < k_t$ then
$P_t^{-1} g$ takes the form $\Phi(a' \cdot x + b')$.
\end{lemma}

Now, Theorem \ref{thmgrad-bound} implies that $\|\E\grad h_t\| \le k_t$
and so by Lemma \ref{lempullback}, $g_t$ is in the range of $P_t$,
and $P_t^{-1} g_t$ is
either the indicator of a half-space or $\Phi$ composed with a linear function.
Let $g = P_t^{-1} g_t$. Then Proposition \ref{proplarge-t} implies that
\[
\E\bigl(P_t (f - g)\bigr)^2 = \E(f_t -
g_t)^2 \le C \frac{\delta(f)^2}{m(f)^C}.
\]
In order to prove Theorem \ref{thmmain}, it suffices
to show that $\E(f - g)^2$ is small. In other words, setting $h = f - g$,
we want to bound $\E h^2$ in terms of $\E(P_t h)^2$. For a general
function $h$,
this is an impossible task. To see why, consider $h_k(x) = \sgn(\sin(kx))$.
Then $\E h_k^2 = 1$ for all $k$, but for any $t > 0$, $P_t h_k \to0$ as
$k \to\infty$. Hence, $\E(P_t h_k)^2 \to0$, and so
$\E h_k^2$ cannot be bounded in terms of $\E(P_t h_k)^2$.

The key to bounding $\E(f - g)^2$ in terms of $\E(P_t (f-g))^2$ is to exploit
some extra information that we have on $f - g$. In particular, we have assumed
that $f$ almost minimizes Bobkov's functional
$\E\sqrt{I^2(f) + \|\grad f\|^2}$. In particular,
\[
\E\|\grad f\| \le\E\sqrt{I^2(f) + \|\grad f\|^2} \le I(
\E f) + \delta(f).
\]
If we assume that $\delta(f) \le1$ (if not, then Theorem \ref
{thmmain} is
meaningless anyway), then $\E\|\grad f\| \le2$. We will translate
this smoothness condition into a condition on the Hermite spectrum
of $f$, which will
allow us to bound $\E(f - g)^2$ in terms of $\E(P_t (f-g))^2$.

We should remark that for nonnegative functions $h$,
reverse hypercontractive inequalities can be used to bound $\E h^2$ in terms
of $\E(P_t h)^2$. The restriction $h \ge0$ prevents the positive and
negative parts of $h$ from canceling out under $P_t$, rendering
examples like
$h_k(x) = \sgn\sin(kx)$ impossible. For our application, however, we
must consider
functions that take positive and negative values.

Our main tool in this section is an inequality by Ledoux, which gives
a connection between $\E\|\grad f\|$ and the action of $P_t$ on $f$.

%th4.2 #&#
\begin{theorem}[(\cite{ledoux94})]\label{thmledoux}
There is a universal constant $C$ such that
for any smooth function $f\dvtx\R^n \to[-1, 1]$
and any $t > 0$,
\[
\E f (f - P_t f) \le C \sqrt{t} \E\|\grad f\|.
\]
\end{theorem}

This inequality (in a sharper form) was originally derived to
show the connection between Borell's noise sensitivity inequality \cite{bor85}
and the Gaussian isoperimetric inequality.
We will give another application:
Theorem \ref{thmledoux} implies that for smooth
functions, the Hermite coefficients decay at a certain rate.

%s4.1 #&#
\subsection{Smoothness and the Hermite expansion}

Recall that the Hermite polynomials
$\{H_\alpha\dvtx \alpha\in\{0, 1, \dots\}^n\}$
form an orthogonal basis of $(\R^n, \gamma_n)$ \cite{Szego75}.
Let
\[
G_\alpha= \frac{H_\alpha}{\sqrt{\E H_\alpha^2}}
\]
be the corresponding
orthonormal basis. We will use the well-known fact
that $P_t$ acts diagonally on this basis:
%
%e25 #&#
\begin{equation}
\label{eqhermite} P_t G_\alpha= e^{-|\alpha|t}
G_\alpha.
\end{equation}

%le4.3 #&#
\begin{lemma}\label{lemhigh-weights}
Suppose $f\dvtx\R^n \to[-1, 1]$ is a smooth function
and $f = \sum_\alpha b_\alpha G_\alpha$. Then
for any $N \in\{1, 2, \dots\}$,
\[
\sum_{|\alpha| \ge N} b_\alpha^2 \le C
N^{-1/2} \E\|\grad f\|.
\]
\end{lemma}

\begin{pf}
By \eqref{eqhermite} and Theorem \ref{thmledoux},
\[
\sum_{\alpha} \bigl(1 - e^{-|\alpha|t}\bigr)
b_\alpha^2 = \E f(f - P_t f) \le C \sqrt t \E\|
\grad f\|.
\]
If $|\alpha| \ge1/t$ then $e^{-|\alpha| t} \le1/e$; hence,
\[
(1 - 1/e)\sum_{|\alpha| \ge1/t} b_\alpha^2
\le\sum_{\alpha} \bigl(1 - e^{-|\alpha|t}\bigr)
b_\alpha^2 \le C \sqrt t \E\|\grad f\|.
\]
Now set $t = \frac{1}{N}$.
\end{pf}

Since we know how the semigroup $P_t$ acts on the Hermite basis
and we know how the Hermite coefficients of nice functions
are distributed, we are in a position to bound
$\E f^2$ in terms of $\E(P_t f)^2$. Essentially,
Lemma \ref{lemhigh-weights} tells us that the high coefficients
do not contribute much to $\E f^2$, while \eqref{eqhermite}
implies that the low coefficients contributing to $\E f^2$
also contribute to $\E(P_t f)^2$.

%le4.4 #&#
\begin{lemma}\label{lemspectral-pullback}
For any smooth $h\dvtx\R^n \to[-1, 1]$ and any $t \ge1$,
\[
\E h^2 \le C \bigl(1 + \E\|\grad h\|\bigr) \sqrt{\frac{t}{\log(1/\E(P_t h)^2)}}.
\]
\end{lemma}

\begin{pf}
Expand $h = \sum_\alpha b_\alpha G_\alpha$
and let $\varepsilon= \E(P_t h)^2$.
Then \eqref{eqhermite} implies that
\[
\varepsilon= \E(P_t h)^2 = \sum
_\alpha e^{-2t|\alpha|} b_\alpha^2.
\]
On the other hand, Lemma \ref{lemhigh-weights} implies that
%
%e26 #&#
\begin{eqnarray}
\label{eqlem-pull-back} \E h^2 &=& \sum_\alpha
b_\alpha^2
\nonumber
\\
&\le& e^{2t(N-1)} \sum_{|\alpha| \le N - 1}
b_\alpha^2 e^{-2t |\alpha|} + \sum
_{|\alpha| \ge N} b_\alpha^2
\\
&\le& e^{2t(N-1)} \varepsilon+ C N^{-1/2} K,\nonumber
\end{eqnarray}
where $K = \E\|\grad h\|$.

Now we choose $N$ to optimize \eqref{eqlem-pull-back}.
Let $\beta= \frac{1}{2t} \log\frac{1}{\varepsilon}$
and set $N = \lceil\beta- \frac{1}{4t} \log\beta\rceil$. Since
$\beta> \log\beta$ and $t \ge1$, $N \ge\beta/2$
(and in particular, $N$ is a positive integer).
Moreover, $N - 1 \le\beta- \frac{1}{4t} \log\beta$ and so (since
$e^{2t\beta} = 1/\varepsilon$)
$e^{2t(N-1)} \varepsilon\le\beta^{-1/2}$.
Plugging these bounds on $N$ back into \eqref{eqlem-pull-back} yields
\[
\E h^2 \le\beta^{-1/2} + C K \beta^{-1/2} \le C(1 +
K) \sqrt{\frac
{t}{\log(1/\varepsilon)}}.
\]
\upqed\end{pf}

%s4.2 #&#
\subsection{Proof of Theorem \texorpdfstring{\protect\ref{thmmain}}{1.3}}
Finally, we are ready to prove Theorem \ref{thmmain}.
As we discussed at the beginning of the section, we may assume that
$\delta= \delta(f) \le1$, which implies that $\E\|\grad f\| \le2$.
We may also assume that $m(f) \ge\log^{-1/2} (1/\delta)$:
if not, then either $\E f \le\log^{-1/2} (1/\delta)$ or
$(1 - \E f) \le\log^{-1/2} (1/\delta)$. In the first case,
$f$ may be approximated well by the zero function, which in turn
may be approximated by functions of the form $\Phi(a \cdot x + b)$.
Specifically, for any $a \in\R^n$,
$\Phi(a\cdot x + b) \to0$ as $b \to-\infty$ and so
\[
\lim_{b \to-\infty} \E\bigl(f(X) - \Phi(a \cdot X + b)
\bigr)^2 = \E f^2 \le \E f \le\frac{1}{\sqrt{\log(1/\delta)}}.
\]
That is, if $\E f \le\log^{-1/2} (1/\delta)$ then the conclusion
of Theorem \ref{thmmain} holds trivially. A similar argument
(but with the zero function replaced by the constant function 1)
holds when $(1 - \E f) \le\log^{-1/2} (1/\delta)$.
Thus, we may assume that $m(f) \ge\log^{-1/2} (1/\delta)$.

As in the discussion at the beginning of the section, take
(by Proposition \ref{proplarge-t}) $t \in[C, C + 1]$ so that
\[
\E(h_t - \E h_t - X \cdot\E\grad h_t)^2
\le C \frac{\delta^{1/4}(f)}{m^{5/4}(f)}.
\]
Let $g_t(x) = \Phi(x \cdot\E\grad h_t + \E h_t)$ and
$g = P_t^{-1} g_t$ (which exists, recall, by Lemma~\ref{lempullback}
and because $\|\E\grad h_t\| \le k_t$).

By Proposition \ref{proplarge-t} and because $\Phi$ is a contraction,
\[
\E(g_t - f_t)^2 \le\E(h_t -
\E h_t - X \cdot\E\grad h_t)^2 \le C
\frac{\delta(f)^{1/4}}{m(f)^{5/4}} \le C\delta^{1/8},
\]
where the last inequality follows from because we have assumed that
$m(f) \ge\log^{-1/2}(1/\delta)$.
Set $h = g - f$. Since $g$ is
$\Phi$ composed with a linear function, $\E\|\grad g\| \le\phi(0)
\le1$,
and hence $\E\|\grad h\| \le\E\|\grad g\| + \E\|\grad f\| \le3$.
By Lemma \ref{lemspectral-pullback},
\[
\E(g - f)^2 \le\frac{C}{\sqrt{\log(1/\E(g_t - f_t)^2)}} \le\frac{C}{\sqrt{\log(1/\delta)}}.
\]
This completes the proof of Theorem \ref{thmmain}.

%s5 #&#
\section{Robust results for sets}\label{secset}

There are two pieces needed to get from Theorem~\ref{thmmain}
to Theorem~\ref{thmmain2}. First, we need to interpret
Theorem \ref{thmmain} in the case that $f$ is an indicator
function (which is necessarily nonsmooth). For this, we
simply apply Theorem \ref{thmmain} to $P_t f$, which is smooth,
and take $t$ to zero. Fortunately for us, most of the work in this
step was done in \cite{ck01}.

%le5.1 #&#
\begin{lemma}[(\cite{ck01})]\label{lemsa}
For any measurable set $A$,
\[
\gamma_n^+(A) = \lim_{t \to0} \E\sqrt{I^2(P_t
1_A) + \|\grad P_t 1_A\|^2}.
\]
\end{lemma}

The second piece we require is something that will let us pass
from a function $\Phi(a \cdot x + b)$ to an affine half-space.
For this piece, we
just round $\Phi(a \cdot x + b)$ to $\{0, 1\}$.

%le5.2 #&#
\begin{lemma}\label{lemto-set}
Let $A$ be a measurable set and
$g(x) = \Phi(a \cdot x + b)$. There exists an affine half-space
$B$ such that
\[
\gamma_n(A \symdiff B) \le\E(1_A - g)^2.
\]
\end{lemma}

\begin{pf}
Let $B = \{x \in\R^n\dvtx a \cdot x + b \ge0\}$.
Since $1_B$ is obtained by rounding $g$ to $\{0, 1\}$, it follows
that $|1_B - g| \le|1_A - g|$ pointwise. Thus,
\[
4 \gamma_n(A \symdiff B) = \E(1_A -
1_B)^2 \le2 \E(1_A - g)^2 + 2
\E(1_B - g)^2 \le4 \E(1_A -
g)^2.\quad%\qedhere
\]
\upqed\end{pf}

\begin{pf*}{Proof of Theorem \ref{thmmain2}}
Let $f = 1_A$ and write $f_t = P_t f$; recall that
$\delta= \gamma_n^+(A) - I(\gamma_n(A))$.
Note that $\E f_t = \E f = \gamma_n(A)$ for all $t > 0$.
Since $\E\sqrt{I^2(f_t) + \|\grad f_t\|^2} \to\gamma_n^+(A)$
by Lemma \ref{lemsa},
\[
\E\sqrt{I^2(f_t) + \|\grad f_t
\|^2} - I(\E f_t) \le2\delta
\]
for all small enough $t$. Now set
$
\varepsilon= \log^{-1/2} (1/\delta)$.
Since the semigroup $P_t$ is strongly continuous in $L_2(\gamma_n)$,
we can take $t$ small enough so that
$\E(f_t - 1_A)^2 \le\varepsilon$.

Apply Theorem~\ref{thmmain} to $f_t$: we receive a function
$g(x) = \Phi(a \cdot x + b)$ with $\E(f_t - g)^2 \le C \varepsilon$.
By the triangle inequality, $\E(1_A - g)^2 \le C \varepsilon$ and
so Lemma \ref{lemto-set} gives us an affine half-space $B$
with $\gamma_n(A \symdiff B) \le C \varepsilon$.
\end{pf*}

%s6 #&#
\section{Conclusion}
%s6.1 #&#
\subsection{Open problems}
To conclude, we present two natural open problems:
\begin{itemize}
\item
Is there a result which strengthens both our result and \cite
{italian11}? For sets $A$ of measure $1/2$, such result should give
the existence of a half space $B$ with $\gamma_n(A \Delta B) \leq C
\delta^{1/2}$ where $C$ is an absolute constant.
\item
Could similar results be obtained for other measures? In particular,
log-concave measures?
This question was suggested to us independently by Franck Barthe,
Michel Ledoux and Shahar Mendelson.
We note that the one-dimensional analogue of \cite{italian11} was
established by \cite{Castro10}.
\end{itemize}

%s6.2 #&#
\subsection{Subsequent work}
Building on techniques that we develop here, we very recently obtained
a robust dimension-free
version of Borell's theorem \cite{MosselNeeman13b}, thereby
establishing the applications
that we mentioned earlier in the \hyperref[sec1]{Introduction}.
Moreover, in Borell's
theorem, we obtained polynomial
rather than logarithmic rates (although we could not achieve the optimal
exponent of $\frac{1}2$). However, our robust version of Borell's theorem
does not imply Theorem \ref{thmmain2} even though Borell's result implies
the isoperimetric inequality as $\rho\to1$, since we were only able
to obtain the correct dependence on $\rho$ for sets $A$ satisfying
$\gamma_n(A) = \frac{1}2$.
%Despite this technical obstacle, our recent
%work \cite{MosselNeeman12b} suggests that our methods can be improved
%to yield a
%polynomial version of Theorem \ref{thmmain2}.

\section*{Acknowledgements}
The authors would like to thank the following colleagues:
Franck Barthe and Almut Buchard for interesting discussions at the
beginning of the project and for comments on
the draft of the manuscript, Michel Ledoux for pointing out several
places where our arguments in the previous Arxived version could be
simplified substantially and for a number of helpful
references and Francesco Maggi for helpful discussions regarding \cite
{italian11}.
We would also like to thank the anonymous referee for offering detailed
and constructive comments.
%\bibliographystyle{plain}
%\bibliography{my,robustness,all-robustness}

% imsref loaded by akundreckaite, 2014-01-15 13:58:37
%

%\begin{appendix}
%\section{}
%\end{appendix}

%\section*{Acknowledgments}

%\begin{supplement}[id=suppA]
%\slink[doi]{10.1214/00-AOPXXXXSUPP} %[doi,text={...}] - jei reikia
%suskaldyti doi
%\sfilename{aopXXXX\_supp.pdf}
%\end{supplement}

%\begin{thebibliography}{99}
%\bibitem[\protect\citeauthoryear{}{}]{r1}
%\bibitem{r1}
%\end{thebibliography}

\printaddresses

\end{document}